\documentclass[11pt]{amsart}

\usepackage[english]{babel}
\usepackage{amsmath}
\usepackage{amsfonts}
\usepackage{amssymb}
\usepackage{amsthm}
\usepackage{epsfig}
\usepackage{graphicx}

\newtheorem{theorem}{Theorem}[section]
\newtheorem{lemma}{Lemma}[section]
\newtheorem{corollary}{Corollary}[section]
\newtheorem{remark}{Remark}[section]
\newtheorem{proposition}{Proposition}[section]
\newtheorem{conjecture}{Conjecture}[section]

\def\conv{\mathop\mathrm{conv}\nolimits}

\def\K{\mathcal{K}}

\def\R{\mathbb{R}}
\def\C{\mathbb{C}}
\def\Real{\mathop\mathrm{Re}\nolimits}
\def\V{\mathrm{V}}

\def\W{\mathrm{W}}
\def\A{\mathrm{A}}
\def\p{\mathrm{p}}
\def\cir{\mathrm{R}}
\def\inr{\mathrm{r}}

\numberwithin{equation}{section}

\begin{document}

\title{Notes on the roots of Steiner polynomials}
\dedicatory{Dedicated to J\"org M.~Wills on the occasion of his 70th birthday}
\author{Martin Henk}
\address{Institut f\"ur Algebra und Geometrie, Otto-von-Guericke
Universit\"at Mag\-deburg, Universit\"atsplatz 2, D-39106-Magdeburg,
Germany} \email{henk@math.uni-magdeburg.de}

\author{Mar\'\i a A. Hern\'andez Cifre}
\address{Departamento de Matem\'aticas, Universidad de Murcia, Campus de
Espinar\-do, 30100-Murcia, Spain} \email{mhcifre@um.es}

\thanks{Second author is supported in part by Direcci\'on General de
Investigaci\'on (MEC) MTM2004-04934-C04-02 and by Fundaci\'on
S\'eneca (C.A.R.M.) 00625/PI/04.}



\begin{abstract}
We study the location and the size of the roots of Steiner
polynomials of convex bodies in the Minkowski relative geometry.
Based on a problem of Teissier on the intersection numbers of Cartier divisors of compact algebraic
varieties it was conjectured that these roots have certain geometric
properties related to the in- and circumradius of the convex body.
We show that the roots of 1-tangential bodies fulfill the
conjecture, but we also present convex bodies violating
each of the conjectured properties.
\end{abstract}

\maketitle

\section{Introduction}

Let $\K^n$ be the set of all convex bodies, i.e., compact convex
sets, in the $n$-dimensional Euclidean  space $\R^n$, and let $B_n$
be the $n$-dimensional unit ball. The subset of $\K^n$ consisting of
all convex bodies with non-empty interior is denoted by $\K_0^n$.
The volume of a set $M\subset\R^n$, i.e., its $n$-dimensional
Lebesgue measure, is denoted by $\V(M)$. For two convex bodies
$K,E\in\K^n$ and a non-negative real number $\rho$ the volume of
$K+\rho\,E$ is a polynomial of degree $n$ in $\rho$ and it can be
written as
\begin{equation}
  \V(K+\rho E) = \sum_{i=0}^n \binom{n}{i} \W_i(K;E)\,\rho^i.
\label{eq:steiner-minkowski}
\end{equation}
This polynomial is called the {\em Minkowski-Steiner polynomial} or
the {\em relative Steiner polynomial} of $K$. The coefficients
$\W_i(K;E)$ are called the {\em relative quermassintegrals} of $K$,
and they are just a special case of the more general defined {\em
mixed volumes} for which we refer to \cite[s.~5.1]{Sch}. In
particular, we have $\W_0(K;E)= \V(K)$, $\W_n(K;E)=\V(E)$ and
$\W_i(K;E)=\W_{n-i}(E;K)$.

If $E=B_n$ the polynomial \eqref{eq:steiner-minkowski} becomes the
classical {\em Steiner polynomial} or {\em Steiner formula}
\cite{Ste40}, and $\W_i(K;B_n)$, for short denoted by $\W_i(K)$, is
the classical {\em $i$-th quermassintegral} of $K$. In this case,
$n\,\W_1(K)$ is the {\em surface area} of $K$,
$\W_n(K)=\V(B_n)=\kappa_n$ is the $n$-dimensional volume of $B_n$
and $(2/\kappa_n)\W_{n-1}(K)$ is the {\em mean width} of $K$
\cite[p.~42]{Sch}.

The {\em relative inradius} $\inr(K;E)$  and {\em relative
circumradius} $\cir(K;E)$ of $K$ with respect to $E$ are defined,
respectively, by
\begin{equation*}
\begin{split}
 \inr(K;E) &= \max\{r : \exists\, x\in\R^n\text{ with } x+r\,E\subseteq K\}, \\
 \cir(K;E) & = \min\{R : \exists\, x\in\R^n\text{ with } K\subseteq x+R\,E\}.
\end{split}
\end{equation*}
Notice that it always holds
\begin{equation}\label{e:inr=1/cir}
\inr(K;E)\, \cir(E;K)=1.
\end{equation}

In the planar case the inradius, circumradius and the
quermassintegrals are related by the well-known Bonnesen inequality
\begin{equation}\label{eq:blascke-inequality}
\W_1(K;E)^2\!-\!\W_0(K;E)\W_2(K;E)\!\geq\!\dfrac{\W_2(K;E)^2}{4}\Bigl(\cir(K;E)-\!\inr(K;E)\Bigr)^2\!.
\end{equation}
Bonnesen \cite{Bo} proved this result for $E=B_2$, the proof of the
general case is due to  Blaschke \cite[pp.~33--36]{Bl49}.
This inequality sharpens (in the plane) the Aleksandrov-Fenchel and
the isoperimetric inequalities and there is no known generalization
of it to higher dimensions. In fact \eqref{eq:blascke-inequality} is
an immediate consequence of the following stronger relations
\cite[pp.~33--36]{Bl49} (see also \cite{Fl})
\begin{equation}\label{eq:roots}
\W_0(K;E)+2\W_1(K;E)\rho+\W_2(K;E)\rho^2\leq 0\; \text{ if }\;
-\cir(K;E)\leq\rho\leq -\inr(K;E).
\end{equation}
The left hand side is just the relative Steiner polynomial
\eqref{eq:steiner-minkowski}, and so \eqref{eq:roots} says that the
(relative) Steiner polynomial
\begin{equation*}
f(K,E,s) = \sum_{i=0}^n \binom{n}{i}\W_i(K;E)\,s^i,
\end{equation*}
regarded as a formal polynomial in a complex variable $s\in\C$, has
in the case $n=2$ two real (negative) roots, one root less than or
equal   $-\cir(K;E)$ and the other root not less than $-\inr(K;E)$.

In \cite{Te82} Teissier studied Bonnesen-Type inequalities in
Algebraic Geometry, more precisely, intersection numbers of Cartier
divisors of $n$-dimensional compact algebraic varieties. These
intersection numbers ``behave similarly'' as the (relative)
quermassintegrals and, in particular, in the case $n=2$ they satisfy
an inequality as \eqref{eq:roots} with suitable defined in- and
circumradius. Teissier raised the problem to find extensions of
these two dimensional properties to higher dimensions (see also
\cite[p.\;103]{Od}). In view of the properties derived from
\eqref{eq:roots} in the planar case, in \cite{SY88n} and
\cite[p.\;65]{SY93} the following conjecture was posed which we
formulate in terms of the relative Steiner polynomial:
\begin{conjecture}\label{co:roots_S-Y}
Let $K,E\in\K^n$. If $a_1\leq\dots\leq a_n$ are the real parts of the
roots of $f(K,E,s)$, then
\[
a_1\leq -\cir(K;E)\leq -\inr(K;E)\leq a_n\leq 0.
\]
\end{conjecture}
As mentioned before, \eqref{eq:roots} implies the conjecture in
dimension 2. A first systematic study of the roots of the classical
Steiner polynomial in the $3$-dimensional case as well as their
relations to the so called Blaschke diagram can be found in
\cite{HCS}. There it is also shown that the conjecture is correct
for some special 3-dimensional convex bodies.  But no further
progress has been made on it; only some results have been obtained
by Wills \cite{W90b} for a closely related polynomial, the so called
{\em Wills functional}.

In the following we will call the conjectured property that all real
parts of the roots of $f(K,E,s)$ are non-positive the {\em negativity property} of
the roots, and the conjectured bounds related to $\inr(K;E)$ and
$\cir(K;E)$ will be referred to as the {\em inradius} and {\em
circumradius bound}, respectively.

In Section 3 we study the above conjecture for the class of
$p$-tangential bodies; for a definition see also Section 3. Among
others we show that $n$-dimensional $1$-tangential bodies, the so
called cap-bodies, verify the conjecture.

\begin{theorem}\label{thm:cap_bodies}
Let $K\in\K^n_0$ be a $1$-tangential body of $E\in\K^n_0$. Then the
roots of $f(K,E,s)$ satisfy Conjecture \ref{co:roots_S-Y}.
\end{theorem}

If we move, however, to $2$-tangential bodies then we loose, in
general, the negativity property of the real parts of the roots.
\begin{theorem}\label{thm:2_tangential}
There exists a $2$-tangential body $K\!\!\in\!\K^{15}\!$ such that
$f(K\!,\!B_{15},s)$ has a root with positive real part.
\end{theorem}
We remark that the minimum dimension of a convex body $K$ such that
$f(K,E,s)$ violates the negativity property of the roots is at least
6. Based on the well-known inequalities
\begin{equation}\label{eq:special_af}
\W_i(K;E)^2 - \W_{i-1}(K;E)\,\W_{i+1}(K;E) \geq 0, \quad 1\leq i\leq
n-1,
\end{equation}
which are particular cases of the Aleksandrov-Fenchel inequality
(see e.g.~\cite[s.~6.3]{Sch}), and on the Routh-Hurwitz criterion
(see e.g. \cite[p.\;181]{Mard}) one can check that $f(K,E,s)$ is a
{\em Hurwitz} polynomial for $n\leq 5$, i.e., all its roots lie in
the left half plane (see also \cite[p.~103]{Te82}).

In Section 4 we construct a $3$-dimensional convex body violating
the conjectured  circumradius bound, more precisely,
\begin{theorem}\label{t:circumradius}
There exists $K\in\K^3_0$ such that all the real parts of the roots
of $f(K,B_3,s)$ are greater than $-\cir(K;B_3)$.
\end{theorem}
Since $\W_i(K,E)=\W_{n-i}(E,K)$ we have $f(K,E,s)=s^n\,f(E,K,1/s)$. Hence,  and on account of
\eqref{e:inr=1/cir} it is not surprising that the body $K$ of the
theorem above leads also to a counterexample for the inradius bound.
\begin{corollary}\label{c:inradius}
There exists $K\in\K^3_0$ such that all the real parts of the roots
of $f(B_3,K,s)$ are less than $-\inr(B_3;K)$.
\end{corollary}

Before giving the proofs of the theorems above we study in Section 2
the size of the roots of the relative Steiner polynomial $f(K,E,s)$.
We prove upper and lower bounds for them in terms of the
circumradius and the inradius.

\section{Bounds for the roots of the Steiner polynomial}
In the following we will use the inequalities
\begin{equation}\label{eq:simple}
\inr(K;E)\,\W_{i+1}(K;E) \leq \W_i(K;E)\leq\cir(K;E)\,\W_{i+1}(K;E),
\end{equation}
for $i\in\{0,\dots,n-1\}$. Since, up to translations,
$\inr(K;E)E\subseteq K$ and $K\subseteq\cir(K;E)E$ these inequalities
are a direct consequence of the monotonicity of the mixed
volumes (cf.~e.g.~\cite[p.~277]{Sch}).  Furthermore, for
a complex number $s$ we denote by $\Real(s)$ its real part.

\begin{proposition}
Let $K\in\K^n$, $E\in\K_0^n$, and let $\gamma_i$, $i=1,\dots,n$, be the roots of
the Steiner polynomial $f(K,E,s)$.
\begin{enumerate}
\item[{\rm i)}] If $\dim K=m$, $m\geq 1$, the non-zero
roots $\gamma_i$ are bounded by
\begin{equation*}
\dfrac{n-m+1}{m}\dfrac{\W_{n-m}(K;E)}{\W_{n-m+1}(K;E)}\leq
|\gamma_i|\leq n\dfrac{\W_{n-1}(K;E)}{\W_n(K;E)}.
\end{equation*}
The upper bound is best possible. In particular, we have
$\inr(K;E)/n\leq|\gamma_i|\leq n\cir(K;E)$.
\item[{\rm ii)}] $\bigl|\Real(\gamma_1)\bigr|+\dots+\bigl|\Real(\gamma_n)\bigr|\geq
n\inr(K;E)$.
\item[{\rm iii)}] If $\Real(\gamma_i)\leq 0$ for all $i=1,\dots,n$,
then
$\bigl|\Real(\gamma_1)\bigr|+\dots+\bigl|\Real(\gamma_n)\bigr|\leq
n\cir(K;E)$.
\end{enumerate}
\end{proposition}

\begin{proof}
Since $\dim K=m$ we have  $\W_i(K;E)=0$ for  $i=0,\dots,n-m-1$ and
$\W_i(K;E)> 0$ for $i=n-m,\dots,n$ (see \cite[p.~277]{Sch}). Hence the Steiner polynomial of $K$
is
\[
f(K,E,s)=s^{n-m}\sum_{i=0}^m\binom{n}{n-m+i}\W_{n-m+i}(K;E)s^i=s^{n-m}h_m(K,E,s)
\]
and the non-zero roots of $f(K,E,s)$ are the roots of $h_m(K,E,s)$.
It is known that the roots of a polynomial
$f(s)=a_0+a_1\,s+\dots+a_m\,s^m$ with positive real coefficients
$a_j$ lie in the ring $\rho_1\leq|s|\leq\rho_2$, where
$\rho_1=\min\{a_j/a_{j+1}\}$ and $\rho_2=\max\{a_j/a_{j+1}\}$, for
$j=0,1,\dots,m-1$, see e.g.~\cite[p.\;137]{Mard}. Hence in the case
of $h_m(K,E,s)$ we just have to find the minimum and maximum of
$\binom{n}{j}\W_j(K;E)/\bigl(\binom{n}{j+1}\W_{j+1}(K;E)\bigr)$,
$j=n-m,\dots,n-1$. By \eqref{eq:special_af} we see that
$\W_j(K;E)/\W_{j+1}(K;E)$ is increasing in $j$, and since
$\binom{n}{j}/\binom{n}{j+1}$ is also increasing we get
\begin{equation*}
\dfrac{n-m+1}{m}\dfrac{\W_{n-m}(K;E)}{\W_{n-m+1}(K;E)}\leq
\dfrac{\binom{n}{j}\W_j(K;E)}{\binom{n}{j+1}\W_{j+1}(K;E)}\leq
n\dfrac{\W_{n-1}(K;E)}{\W_n(K;E)},
\end{equation*}
for $j=n-m,\dots,n-1$, which shows the inequalities in i). The bounds in
i) in terms of the inradius and the circumradius follow immediately with
\eqref{eq:simple}. Notice that if $m<n$ then $\inr(K;E)=0$. Hence the only
non-trivial lower bound is obtained when $m=n$.

By i) we see that the only non-zero root of the Steiner polynomial of a
line segment $K$, i.e., $\dim K=1$,  is given by
$\gamma=-n\W_{n-1}(K;E)/\W_n(K;E)$. This shows that the upper bound in i)
is best possible for any choice of the body $E\in\K_0^n$.

Since $f(K,E,s)=\W_n(K;E)\prod_{i=1}^n(s-\gamma_i)$ we find that
\begin{equation*}
\gamma_1+\dots+\gamma_n=-n\dfrac{\W_{n-1}(K;E)}{\W_n(K;E)}.
\end{equation*}
Together with \eqref{eq:simple} we get
\begin{equation*}
n\inr(K;E)\leq
|\gamma_1+\dots+\gamma_n|=\bigl|\Real(\gamma_1)+\dots+\Real(\gamma_n)\bigr|\leq
\bigl|\Real(\gamma_1)\bigr|+\dots+\bigl|\Real(\gamma_n)\bigr|,
\end{equation*}
which shows ii). Finally, if all the roots have negative real part (as conjectured),
we even have
$\bigl|\Real(\gamma_1)\bigr|+\dots+\bigl|\Real(\gamma_n)\bigr|
=\bigl|\Real(\gamma_1)+\dots+\Real(\gamma_n)\bigr|$ and as above we
get by \eqref{eq:simple} the  upper bound stated in iii).
\end{proof}

\section{Roots of the Steiner polynomial and tangential bodies}\label{s:Teissier}

A convex body $K\in\K^n$ containing the convex body $E\in\K^n$ is
called a $p$-tangential body of $E$, $p\in\{0,\dots,n-1\}$, if each
support plane of $K$ that is not a support plane of $E$ contains
only ($p-1$)-singular points of $K$ \cite[p.\;76]{Sch}. Here a
boundary point $x$ of $K$  is said to be an $r$-singular point of
$K$ if the dimension of the normal cone in $x$ is at least $n-r$.
For further characterizations and properties of $p$-tangential
bodies we refer to \cite[Section 2.2]{Sch}.

So a $0$-tangential body of $E$ is just the body $E$ itself and each
$p$-tangential body of $E$ is also a $q$-tangential body for
$p<q\leq n-1$. A $1$-tangential body is usually called {\it
cap-body}, and it can be seen as the convex hull of $E$ and
countably many points such that the line segment joining any pair of
those points intersects $E$.

If $K$ is a $p$-tangential body of $E$ then $\inr(K;E)=1$, and the
following theorem gives a characterization of $n$-dimensional
$p$-tangential bodies in terms of the quermassintegrals
(cf.~\eqref{eq:simple}).
\begin{theorem}[Favard \cite{Fav33}, {\cite[p. 367]{Sch}}]\label{t:Favard}
Let $K,E\in\K^n_0$ and $p\in \{0,\dots,n-1\}$. Then
$\W_0(K;E)=\W_1(K;E)=\cdots=\W_{n-p}(K;E)$ if and only if $K$ is a
$p$-tan\-gen\-tial body of $E$.
\end{theorem}

The next proposition shows that $n$-dimensional $p$-tangential
bodies always satisfy the inradius bound of  Conjecture
\ref{co:roots_S-Y}.
\begin{proposition}\label{p:inradius_bound}
Let $p\in\{0,\dots,n-1\}$, and let $K\in\K^n_0$ be a $p$-tangential
body of $E\in\K^n_0$. Then there  exists a root $\gamma$ of
$f(K,E,s)$ such that $\Real(\gamma)\geq-\inr(K;E)$.
\end{proposition}

\begin{proof}
Let $\gamma_i$, for $1\leq i\leq n$,  be the roots of $f(K,E,s)$.
From $f(K,E,s)=\W_n(K;E)\prod_{i=1}^n(s-\gamma_i)$ we get
\begin{equation*}
\begin{split}
(-1)^n\W_0(K;E) & =\W_n(K;E)\,\gamma_1\cdot\ldots\cdot\gamma_n\;\;\text{ and }\\
(-1)^{n-1}n\W_1(K;E) & =\W_n(K;E)\sum_{i=1}^n \prod_{j\ne i} \gamma_j.
\end{split}
\end{equation*}
By Theorem \ref{t:Favard} we have $\W_0(K;E)=\W_1(K;E)$ for any
$p$-tangential body, and thus
\[
-n=-n\dfrac{\W_1(K;E)}{\W_0(K;E)}=\dfrac{1}{\gamma_1}+\dots+\dfrac{1}{\gamma_n}
=\Real\left(\dfrac{1}{\gamma_1}\right)+\dots+\Real\left(\dfrac{1}{\gamma_n}\right).
\]
Therefore, there exists a root $\gamma_j$, say, such that
$\Real(1/\gamma_{j})\leq -1$, and so $\Real(\gamma_{j})\geq
-1=-\inr(K;E)$.
\end{proof}

\begin{remark}
Let $\gamma_j$ with $\Real(1/\gamma_{j})\leq -1$ be the root of the
above proof of the polynomial $f(K,E,s)$. Since
$\W_i(K;E)=\W_{n-i}(E;K)$ and thus $f(K,E,s)=s^n\,f(E,K,1/s)$ we see
that $1/\gamma_{j}$ is a root of $f(E,K,s)$.  By \eqref{e:inr=1/cir}
we have $R(E;K)=1$ and so we get that the polynomial $f(E,K,s)$
satisfies the circumradius bound.
\end{remark}

Now we come to the proof of Theorem \ref{thm:cap_bodies} and show
that  $1$-tangential bodies (cap-bodies) fulfill the conjecture. We
remark  that the analogous $3$-dimensional result for $E=B_3$ was
already mentioned in \cite[p.\;65]{SY93}.

\begin{proof}[Proof of Theorem \ref{thm:cap_bodies}]
If $K\in\K^n_0$ is a cap-body of $E\in\K^n_0$  Theorem
\ref{t:Favard} asserts that $\W_0(K;E)=\W_i(K;E)$, for all
$i=1,\dots,n-1$, and so we can rewrite the Steiner polynomial in the
following way
\[
\begin{split}
f(K,E,s) & =\W_0(K;E)\left[\sum_{i=0}^{n-1}\binom{n}{i}\,s^i+\dfrac{\W_n(K;E)}{\W_0(K;E)}\,s^n\right]\\
 & =\W_0(K;E)\Bigl[(1+s)^n-\bigl(1-\alpha(K,E)\bigr)\,s^n\Bigr],
\end{split}
\]
where $\alpha(K,E)=\W_n(K;E)/\W_0(K;E)$. Observe that
$0<\alpha(K,E)\leq 1$ since $E\subseteq K$. So in this case  all the
roots $\gamma_k$ of $f(K,E,s)$ have to satisfy the equation
$(1/s+1)^n=1-\alpha(K,E)$ and so we obtain
\begin{equation*}
\dfrac{1}{\gamma_k}=-1+\sqrt[n]{1-\alpha(K,E)}\,\,\mathrm{e}^{\frac{2\pi
k}{n}\mathrm{i}}
\end{equation*}
for $k=0,\dots,n-1$. Hence all the real parts of $1/\gamma_k$, and
thus of $\gamma_k$ are non-positive which shows the negativity
property of the roots stated in Conjecture \ref{co:roots_S-Y}.
The inradius bound of the conjecture is guaranteed by Proposition
\ref{p:inradius_bound} and it remains to verify the circumradius
bound. The smallest real part among the roots $\gamma_k$
corresponds to $k=0$, i.e., it is
\[
\gamma_0=\dfrac{-1}{1-\sqrt[n]{1-\alpha(K,E)}},
\]
and we have to show  that
$-1/\bigl(1-\sqrt[n]{1-\alpha(K,E)}\bigr)\leq -\cir(K;E)$, which is
equivalent to
\begin{equation}\label{eq:to_prove}
1-\dfrac{\W_n(K;E)}{\W_0(K;E)}\geq\left(1-\frac{1}{\cir(K;E)}\right)^n.
\end{equation}
In order to prove \eqref{eq:to_prove}  we use the following inequality
\begin{equation*}\label{e:Sangwine-Yager}
\left(\!\dfrac{\W_{n-1}(K;E)}{\W_0(K;E)}\!\right)^{\frac{n}{n-1}}\!\!-\dfrac{\W_n(K;E)}{\W_0(K;E)}\geq\!
\left[\left(\!\dfrac{\W_{n-1}(K;E)}{\W_0(K;E)}\!\right)^{\frac{1}{n-1}}\!\!-\frac{1}{\cir(K;E)}\right]^n,
\end{equation*}
obtained in \cite[Corollary 22]{SY88n}. Since
$\W_{n-1}(K;E)=\W_0(K;E)$ we get immediately \eqref{eq:to_prove}.
\end{proof}

In contrast to $1$-tangential bodies, the roots of Steiner
polynomials of $2$-tangential bodies do not fulfill all the
properties of Conjecture \ref{co:roots_S-Y}, more precisely, some of
their roots can have a positive real part if the dimension is large
enough.

\begin{proof}[Proof of Theorem \ref{thm:2_tangential}]
On account of Theorem \ref{t:Favard} we may write the Steiner
polynomial $f(K,B_n,s)$ of a $2$-tangential body $K\in\K^n_0$ of
$B_n$ as
\begin{equation*}
f(K,B_n,s)=\W_0(K)\left[\sum_{i=0}^{n-2}\binom{n}{i}\,s^i+n\dfrac{\W_{n-1}(K)}{\W_0(K)}s^{n-1}+\dfrac{\W_n(K)}{\W_0(K)}\,s^n\right].
\end{equation*}
We write $\beta(K)=\W_{n-1}(K)/\W_0(K)$ and
$\alpha(K)=\W_n(K)/\W_0(K)$. All the roots of such a polynomial are
non-zero. Since we are only interested in the negativity property of
its roots we may divide by $s^n$ and with $\mu=1/s$ it suffices to
consider the polynomial
\begin{equation*}
h(K,\mu)=\sum_{i=2}^n\binom{n}{i}\mu^i + n\,\beta(K)\mu +\alpha(K).
\end{equation*}
Now it can be checked with a computer or by applying the Routh-Hurwitz
criterion that the polynomial $\sum_{i=2}^n\binom{n}{i}\mu^i$ has a root
with positive real part for $n=15$. Hence, if we find a 2-tan\-gential
body $K\in\K^{15}_0$ for which $\beta(K)$ and $\alpha(K)$ can be
arbitrarily small we get a counterexample to the negativity property of
Conjecture \ref{co:roots_S-Y}. Observe that the roots of a polynomial are
continuous functions of the coefficients of the polynomial
(cf.~e.g.~\cite[p.~3]{Mard}).

In order to construct such a body let $P_\lambda \in\K^3_0$, $\lambda\geq
2$, be the pyramid over a square basis with vertices
\begin{equation*}
(\pm\lambda,\pm\lambda,-1)^\intercal,\quad
\bigl(0,0,1+2/(\lambda^2-1)\bigr)^\intercal.
\end{equation*}
The coordinates are chosen such that the largest ball contained in
$P_\lambda$ is $B_3$ and that all $2$-faces (facets) of $P_\lambda$
touch $B_3$. Next we embed $P_\lambda$ in the canonical way into
$\R^{15}$ and consider
$K_\lambda=\conv\{P_\lambda,B_{15}\}\in\K^{15}_0$. If $H$ is a
support plane of $K_\lambda$ which is not a support plane of
$B_{15}$ it must be a support plane of $P_\lambda$ and it can not
contain any of the $2$-faces of $P_\lambda$. Thus $H$ contains only
$1$-singular points of $K_\lambda$ and this shows that $K_\lambda$
is a $2$-tangential body of $B_{15}$.

It is easy to see that for the pyramid $P_{\lambda}$ there exists a
constant $c$ such that its 3-dimensional volume is not smaller than
$c\,\lambda^2$. Hence there exists a constant $c_n$ depending only on the
dimension such that $\V(K_\lambda)=\W_0(K_\lambda)\geq c_n\,\lambda^2$. On
the other hand, the circumradius of $K_\lambda$ is certainly less than
$2\lambda$ and so by \eqref{eq:simple} we have the bound
$\W_{n-1}(K_\lambda)\leq 2\lambda \V(B_{15})=2\lambda\,\W_n(K_\lambda)$.
Thus
\begin{equation*}
\lim_{\lambda\to\infty} \beta(K_\lambda)= \lim_{\lambda\to\infty} \alpha(K_\lambda)=0,
\end{equation*}
which shows that $h(K_\lambda,\mu)$ and thus $f(K_\lambda,B_{15},s)$
has a root with positive real part if $\lambda$ is large enough.
\end{proof}

\begin{remark}
It can be also checked that the polynomial $\sum_{i=3}^n
\binom{n}{i} \mu^i$ has a root with positive real part in dimension
$n=12$. Hence by an analogous argument and construction  as above
one can show that there exists a $3$-tangential body $K\in\K^{12}_0$
of $B_{12}$ violating the negativity property of Conjecture
\ref{co:roots_S-Y}.
\end{remark}

Since the construction of this 12-dimensional counterexample  is a
bit more involved we omit it here. To close the section we give a
$20$-dimensional numerical counterexample based on an almost regular
crosspolytope.

For $0<\lambda_1\leq\dots\leq\lambda_n$ we denote by
$C_n^*(\lambda_1,\dots,\lambda_n)$ the orthogonal crosspolytope given by
$C_n^*(\lambda_1,\dots,\lambda_n)=\conv\{\pm \lambda_i\,e_i:
i=1,\dots,n\}$, where $e_i$ denotes the $i$-th canonical unit vector.
Analogously to the proof of Lemma 2.1 in \cite{BH93} it can be shown:
\begin{lemma}
Let
$F^i(\lambda_{l_1},\dots,\lambda_{l_{i+1}})=\conv\{\lambda_{l_1}e_{l_1},\dots,\lambda_{l_{i+1}}e_{l_{i+1}}\}$,
$0\leq i\leq n$, be an $i$-dimensional face of
$C_n^*(\lambda_1,\dots,\lambda_n)$, $1\leq l_1<\dots<l_{i+1}\leq n$.
The external angle of the face
$F^i(\lambda_{l_1},\dots,\lambda_{l_{i+1}})$ is given by
\[
\dfrac{2^{n-i-1}}{\pi^{(n-i)/2}}\dfrac{\sqrt{\sum_{k=1}^{i+1}\prod_{j\neq
k}\lambda_{l_j}^2}}{\prod_{j=1}^n\lambda_j}
\!\int_0^{\infty}\!\mathrm{e}^{-\left(\sum_{k=1}^{i+1}\frac{1}{\lambda_{l_k}^2}\right)x^2}\!\!\left(\prod_{j\neq\{l_1,\dots,l_{i+1}\}}
\int_0^{x}\mathrm{e}^{-\frac{y^2}{\lambda_j^2}}dy\right)\,dx.
\]
\end{lemma}

The $i$-face $F^i(\lambda_{l_1},\dots,\lambda_{l_{i+1}})$ is the
$i$-simplex
$\conv\{\lambda_{l_1}e_{l_1},\dots,\lambda_{l_{i+1}}e_{l_{i+1}}\}$,
and its $i$-th volume is given by
\begin{equation*}
\dfrac{1}{i!}\sqrt{\sum_{k=1}^{i+1}\prod_{j\neq k}\lambda_{l_j}^2}.
\end{equation*}
Since $C_n^*(\lambda_1,\dots,\lambda_n)$ has $2^{i+1}$ equal
$i$-faces of the type $F^i(\lambda_{l_1},\dots,\lambda_{l_{i+1}})$
we get the following formulae for the quermassintegrals.
\begin{theorem}\label{t:V_i_cross-polyt}
Let $0<\lambda_1\leq\dots\leq\lambda_n$. The quermassintegrals of the
ortho\-gonal crosspolytope $C_n^*(\lambda_1,\dots,\lambda_n)$ are given
by:

\medskip

$\W_0\bigl(C_n^*(\lambda_1,\dots,\lambda_n)\bigr)=\dfrac{2^n}{n!}\lambda_1\cdot...\cdot\lambda_n$,

$\W_1\bigl(C_n^*(\lambda_1,\dots,\lambda_n)\bigr)=\dfrac{2^n}{n!}\sqrt{\displaystyle\sum_{i=1}^n\prod_{j\neq
i}\lambda_j^2}$,

\medskip

\noindent and for $0\leq i\leq n-2$,
\[
\begin{split}
\W_{n-i}\bigl(C_n^*(\lambda_1,\dots,\lambda_n)\bigr) &
    =\dfrac{\kappa_{n-i}}{\binom{n}{i}}\dfrac{2^n}{i!\pi^{(n-i)/2}}\sum_{1\leq l_1<\dots<l_{i+1}\leq n}
    \left[\vphantom{\sqrt{\sum^i_{\substack{j\\j}}}}
    \dfrac{\sum_{k=1}^{i+1}\prod_{j\neq k}\lambda_{l_j}^2}{\prod_{j=1}^n\lambda_j}\right.\\
 & \;\left.\int_0^{\infty}\mathrm{e}^{-\left(\sum_{k=1}^{i+1}\frac{1}{\lambda_{l_k}^2}\right)x^2}\!\left(\prod_{j\neq\{l_1,\dots,l_{i+1}\}}
    \int_0^x\mathrm{e}^{-\frac{y^2}{\lambda_j^2}}dy\right)dx \vphantom{\sqrt{\sum^i_{\substack{j\\j}}}}\right].
\end{split}
\]
\end{theorem}
Numerical computations show that the Steiner polynomial of the
cross\-polytope $C_n^*\bigl(\lambda,\overset{{\small
(n/2)}}{\dots},\lambda,1,\overset{{\small (n/2)}}{\dots},1\bigr)$,
$n$ even, has roots with positive real parts for different values of
$\lambda$ and $n$, e.g., for $n=20$ and $\lambda=0.01$.

\section{The in- and circumradius bound}

Now we deal with the part of Conjecture \ref{co:roots_S-Y} regarding
the inradius and the circumradius bounds and we are going to show
that the conjectured properties do not hold in dimension 3. To this
end we consider first planar convex bodies in $\R^3$ and show that
the circumradius bound is false in general for that class of bodies.
From such a planar body  we can easily construct a $3$-dimensional
counterexample which also fails the inradius bound. The inradius and
circumradius of a convex body with respect to the unit ball will be
denoted by $\inr(K)$ and $\cir(K)$. For a planar convex body we
denote, as usual, by $\A(K)$ and $\p(K)$ its area and perimeter,
respectively. With this notation \eqref{eq:roots} can be rewritten
in the case $E=B_2$ as
\begin{equation}\label{eq:bonnesen-ineq}
\A(K)+\p(K)\,\rho+\pi\,\rho^2\leq 0,\; \text{ if } -\cir(K)\leq\rho\leq
-\inr(K).
\end{equation}
The next lemma gives a characterization of those planar convex
bodies in $\R^3$ failing the circumradius bound.

\begin{lemma}\label{lem:dim2}
Let $K\in\K^3$ be a planar convex body. All the roots of its Steiner
polynomial $f(K,B_3,s)$ have real part greater than $-\cir(K)$ if and only
if $\p(K)^2<\bigl(128/(3\pi)\bigr)\,\A(K) $ and $\p(K)<(16/3)\,\cir(K)$.
\end{lemma}

\begin{proof}
Since $K\!\in\!\K^3$ is a planar body we have $\W_0(K)\!=\!0$,
$3\W_1(K)=2\A(K)$ and $3\W_2(K)=(\pi/2)\p(K)$. Hence
\begin{equation*}
f(K,B_3,s)=
s\left(2\A(K)+\frac{\pi}{2}\p(K)\,s+\frac{4}{3}\pi\,s^2\right).
\end{equation*}
The non-zero roots $\gamma_{1,2}$ are given by
\begin{equation*}
\gamma_{1,2}=3\frac{-\p(K)\!\pm\!\sqrt{\p(K)^2\!-\!\frac{128}{3\pi}\A(K)}}{16}.
\end{equation*}
If $\p(K)^2-\bigl(128/(3\pi)\bigr)\,\A(K)\geq 0$, i.e., if $f(K,B_3,s)$
has only real roots, then we have
\begin{equation*}
3\frac{-\p(K)\!-\!\sqrt{\p(K)^2\!-\!\frac{128}{3\pi}\A(K)}}{16}\leq\frac{-\p(K)-\sqrt{\p(K)^2-4\pi\A(K)}}{2\pi}
\leq -\cir(K),
\end{equation*}
where the last inequality follows from \eqref{eq:bonnesen-ineq}.
Hence  all the roots of the Steiner polynomial of $K$ have real
parts greater than $-\cir(K)$ if and only if
$\p(K)^2-\bigl(128/(3\pi)\bigr)\,\A(K)<0$ and
$-(3/16)\p(K)>-\cir(K)$.
\end{proof}

\begin{figure}[h]
\begin{center}
\includegraphics[width=4cm]{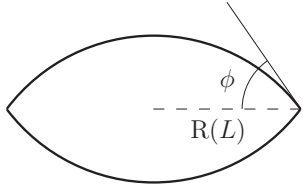}
\caption{A symmetric lens $L$.}\label{f:lens}
\end{center}
\end{figure}

In order to find a convex body satisfying the conditions of the
previous lemma we consider a symmetric lens $L$ with circumradius
$\cir(L)=1$ and perimeter $\p(L)=5.2<16/3$. It is well known (see
e.g.~ \cite[p.\;88--89]{BF}) that symmetric lenses are the extremal
sets of  the inequality
\begin{equation*}
8\phi\,\A(K)\leq\p(K)\bigl(\p(K)-4\cir(K)\cos\phi\bigr),
\end{equation*}
where $\phi$ is determined as the positive solution of the equation
$\p(K)\sin\phi=4\cir(K)\phi$. In the case of a lens, $\phi$ is the half
angle between the two circles of the lens, see Figure \ref{f:lens}.

For our lens $L$ we obtain by the equality case in the above inequality
$\A(L)=2.038627\dots$, which satisfies also the second inequality
$\p(K)^2<\bigl(128/(3\pi)\bigr)\A(K)$ of Lemma \ref{lem:dim2}.

Hence $L$ is a $2$-dimensional convex body in $\R^3$ violating the
circumradius bound. The roots of $f(L,B_3,s)$ are given by
\begin{equation}
\label{eq:roots_lens}
\gamma_{1,2}=-0.975\pm 0.150823\dots \mathrm{i}\quad \text{ and
}\quad \gamma_3=0,
\end{equation}
whose real parts are greater than $-\cir(L)=-1$.

From such a set  it is easy to obtain 3-dimensional convex bodies
$K$ for which all the roots of $f(K,B_3,s)$  have real part greater
than $-\cir(K)$.

\begin{proof}[Proof of Theorem \ref{t:circumradius}]
Observe, in general, we have for $\nu,\,\rho \geq 0$ that $\V(K+\nu\,E+
\rho\,E)=\V\bigl(K+(\nu+\rho)\,E\bigr)$ and thus
\begin{equation*}
f(K+\nu\,E,E,s)=\sum_{i=0}^n \binom{n}{i} \W_i(K;E)\,(s+\nu)^i =
f(K,E,s+\nu).
\end{equation*}
Hence the roots of $f(K+\nu\,E,E,s)$ are given by $\gamma_i-\nu$, where
$\gamma_i$ are the roots of $f(K,E,s)$. Now let $L_\nu=L+\nu B_3$, where
$L$ is the above $2$-dimensional lens. Then $L_\nu\in\K_0^3$ for $\nu>0$
and the roots of $f(L_\nu,B_3,s)$ are $\gamma_i-\nu$, where $\gamma_i$ are
given in \eqref{eq:roots_lens}. Thus all the real parts of $\gamma_i-\nu$
are greater than $-1-\nu=-\cir(L_\nu,B_3)$.
\end{proof}

Finally we deal with the inradius bound. To this end let $L_\nu$ be as in
the proof above. Since $\W_i(L_\nu;B_3)=\W_{3-i}(B_3;L_\nu)$ we get
\begin{equation*}
\begin{split}
f(B_3,L_{\nu},s)  &
    =\left(2\A(L)\nu+\frac{\pi}{2}\p(L)\nu^2+\frac{4}{3}\pi\nu^3\right)s^3\\
  & \quad+\left(2\A(L)+\pi\p(L)\nu+4\pi\nu^2\right)s^2+\left(\frac{\pi}{2}\p(L)+4\pi\nu\right)s +\frac{4}{3}\pi.
\end{split}
\end{equation*}
For instance, for $\nu=1$ it can be checked that the roots of this
polynomial are
\begin{equation*}
\gamma_{1,2}=-0.503393\dots\pm 0.038442\dots \mathrm{i}\quad \text{
and }\quad \gamma_3=-1,
\end{equation*}
whose real parts are smaller than
$-1/2=-1/\cir(L_1,B_3)=-\inr(B_3,L_1)$.

This provides the required  counterexample to the inradius bound of
Conjecture \ref{co:roots_S-Y} stated in Corollary \ref{c:inradius}.

We remark that the above construction in the proof of Theorem
\ref{t:circumradius} and Corollary \ref{c:inradius} works with any
symmetric lens $L$ with $\cir(L)=1$ and perimeter $\p(L)\in[\p_0,16/3]$,
where $\p_0=4\phi_0/\sin\phi_0=5.052\dots$ and $\phi_0$ is the smallest
positive solution of $\phi-\sin\phi\cos\phi=(3\pi/16)\phi^2$.

\medskip

\noindent {\it Acknowledgements.} This work has been developed
during a research stay of the second author at the Otto-von-Guericke
Universit\"at Magdeburg, Germany, supported by Programa Nacional de
ayudas para la movilidad de profesores de universidad e
investigadores espa\~noles y extranjeros (MEC), Ref. PR2006-0351.

\end{document}